\documentclass[10pt]{article}    % Specifies the document style.  Here 
\usepackage{amssymb,latexsym,amsmath} % See Gratzer: "Math Into \LaTeX"
\usepackage[mathscr]{eucal}  % See Gratzer P. 192.
\newfont{\sffl}{msbm10 at 16pt} % See PC TeX User Manual p.112 & p. 117
\newfont{\sff}{msbm10 at 10pt}

\begin{document}           % End of preamble and beginning of text.

\title{Four Properties of Reproducing Kernel Hilbert Spaces\thanks{\small{Approved for public release; distribution is unlimited.}}
}  % Declares the document's title.
 
\author{Alan Rufty%\\         % Declares the author's name.
%\\
%P.O. Box 711\\
%Dahlgren, VA. 22448
}
\date{November 28, 2007}

\maketitle                 % Produces the title.

\newcommand{\KD}{K_{\text{D}}}
\begin{abstract}% See Gratzer: "Math Into \LaTeX", p. 46. Shoot for < 251 words.
  A reproducing kernel Hilbert space (RKHS) has four well-known easily derived properties.  Since these properties are usually not emphasized as a simple means of gaining insight into RKHS structure, they are singled out and proved here.
\end{abstract}

\newcommand{\ls}{\vphantom{\big)}} 
\newcommand{\lsm}{\!\vphantom{\big)}} 
\hyphenation{DIDACKS}

\vskip .05in
\noindent
\begin{itemize}
\item[\ \ ] \small{\textbf{Key words:} {reproducing kernel, Dirichlet form}}
\item[\ \ ] \small{\textbf{AMS subject classification (2000):} {Primary 46E22.}}
\end{itemize}

\renewcommand {\baselinestretch}{1.35} % See Gratzer p. 105
 
\section{Introduction}\label{S:intro}

   A recent article by the author built on the concept of a reproducing kernel Hilbert space (RKHS) \cite{DIDACKS} and the present article provides additional standard background material on RKHSs.  For concreteness, since the primary norms and inner products of interest in \cite{DIDACKS} were related to Dirichlet integrals defined over some connected $\mathbb{R}^3$ region $\Omega$, the same overall setting and notation is assumed here, although the properties and proofs given apply to more general settings; hence, the admissible functions and inner products are assumed to be real valued and vectors denote points of $\mathbb{R}^3$.  While general introductions to RKHS theory can be found in \cite{Aronszajn}, \cite{Moritz} or \cite{Mate}, the present discussion is limited in scope, but is more-or-less self contained and should be accessible to a wide readership.

\section{Reproducing Kernel Properties}\label{S:pmath}

 First consider RKHS theory where, as it will be presently shown, all reproducing kernels are symmetric.  Let points $\vec{P}$ and $\vec{Q}$ be $\mathbb{R}^3$ points in the same connected region ($\Omega$) and let $\mathscr{H}$ denote the associated Hilbert space of real valued functions defined over $\Omega$ with a real valued inner product, $(\,\,\cdot\,\,,\,\,\cdot\,\,)$.  A reproducing kernel $K(\vec{P},\,\vec{Q})$ associated with $\mathscr{H}$ can be compactly characterized by the following two requirements \cite{Moritz}:
\renewcommand{\theenumi}{\Roman{enumi}}  % See The LaTeX Companion by Gossens, Mittelbach Samarin P. 57
\renewcommand{\labelenumi}{(\theenumi)} % See ibid P. 58
\begin{enumerate}
\item
$K(\vec{P},\,\,\cdot\,) \in \mathscr{H}\ $ and $K(\,\cdot\,,\,\vec{Q}) \in \mathscr{H}\,\,$, which is to say $\|K\|$ must be bounded when treated as a function of either argument. 
\item
$\boldsymbol{(}K{(\vec{P},\vec{X})},\,f(\vec{X})\,\boldsymbol{)} \equiv f(\vec{P})$.%$(K(\vec{P},\,\,\cdot\,)\,,\,f) = f(\vec{P})$
\end{enumerate}
As an aside, observe that while a Dirac delta function satisfies the second requirement, it fails to satisfy the first one so it is not a reproducing kernel.
The actual form and existence of a closed form reproducing kernel is closely tied to the shape of the region of interest.  Reproducing kernels and their associated spaces (RKHSs) were first studied by Bergman and others and then brought to a mature state of development over fifty years ago (for a summary of this work and historical comments see Aronszajn \cite{Aronszajn}, or \cite{RKHSbook} \{which contains \cite{Aronszajn}\}).  

  Next consider the following four important properties of a reproducing kernel:
\newcounter{mycount}  %see Hahn (LaTex for Everyone) p 313 (pp 211-213). Also see pp.232-234.
\begin{list}{\textbullet \ \ Property (\arabic{mycount})}% For \textbullet see Gratzer p 85.
{\usecounter{mycount} \itemsep 0in \parsep .1in \labelwidth 2.0in \leftmargin 1.0in \rightmargin 1.0in} %\leftmargin .45in \rightmargin \leftmargin}
\item %[Property (1)]
 $K{(\vec{P},\vec{Q})} \equiv K{(\vec{Q},\vec{P})}$; i.e., all real reproducing kernels are symmetric.
\item %[Property (2)]
 $K(\vec{P},\vec{Q})$ is bounded.
\item %[Property (3)]
 Two different reproducing kernels over $\Omega$ cannot exist for the same norm. 
\item %[Property (4)]
 For norms that can be expressed as integrals over $\Omega$, two different norms cannot share the same reproducing kernel. 
\end{list}

\noindent
Let $\vec{X}$ denote the dummy integration or inner-product variables.  Property (1) follows immediately from the fact that $(\,K(\vec{P},\,\vec{X})\,,K(\vec{Q},\,\vec{X})\,) = K(\vec{Q},\,\vec{P})$\,, $(\,K(\vec{Q},\,\vec{X})\,,K(\vec{P},\,\vec{X})\,) = K(\vec{P},\,\vec{Q})$ and that $(f,\,g) = (g,\,f)$.  Since this property means that every reproducing kernel is symmetric, the adjective symmetric will normally be used.  Next consider Property (2).  Let $\mathscr{H}$ be a space of square integrable functions so that $\|f\|$ is bounded for all $f \in \mathscr{H}$, then $ \|K(\vec{P},\,\,\cdot\,)\,\|^2 :=(\,K(\vec{P},\,\,\cdot\,),\,K(\vec{P},\,\,\cdot\,)\,) = K(\vec{P},\,\vec{P})$ is bounded.  Because $\|K(\vec{P},\,\,\cdot\,)\,\|$ and $\|K(\vec{Q},\,\,\cdot\,)\,\|$ are bounded, it follows that $K(\vec{P},\,\vec{Q})$ is bounded since $K(\vec{P},\,\vec{Q}) = (\,K(\vec{P},\,\,\cdot\,),\,\,K(\vec{Q},\,\,\cdot\,)\,\,) \leqq  \|K(\vec{P},\,\,\cdot\,)\|\,\|K(\vec{Q},\,\,\cdot\,)\|$.  To prove Property~(3) assume to the contrary that some $K_{A}$ and $K_{B}$ exists for a given norm with $K_{A} \neq K_{B}$.  Then $K_{A}{(\vec{P},\vec{Q})} = (K_{A}{(\vec{P},\vec{X})},\,K_{B}{(\vec{X},\vec{Q})}) = K_{B}{(\vec{P},\vec{Q})}$; consequently, Property~(3) must hold.

  To prove Property~(4) some additional notation is required. First consider only inner-products defined without embedded operators that can be described in terms of weighted integrals:  $\boldsymbol{(}f,\,g\boldsymbol{)} = \int fg\,{\mu}(\vec{X})\,\,d^3\,\vec{X}$, where $\mu$ is the weight.  Let two such different norms or inner products exist and label them $\boldsymbol{(}\,\cdot\,,\,\cdot\,\boldsymbol{)}{\lsm}_A$ and $\boldsymbol{(}\,\cdot\,,\,\cdot\,\boldsymbol{)}{\ls}_B$.  The statement that these yield different norms then means that $f$ and $g$ always exist such that $\boldsymbol{(}f,\,g\boldsymbol{)}{\lsm}_A \neq \boldsymbol{(}f,\,g\boldsymbol{)}_B$\,, which, in turn, means that since the regions are the same they have different weight functions:  ${\mu}_A \neq {\mu}_B$.  Ordinarily a vector dummy argument is tied to the implementation of the norm and the choice of symbol for it does not matter, but here the choice of dummy arguments must be tracked, so for $\nu = A$ or $B$ let $\boldsymbol{(}f,\,g\boldsymbol{)}_{\nu(\vec{Q})}$ indicate that vector field variable $\vec{Q}$ fills this role.  Then, contrary to our assumption, observe that
\begin{multline}\notag
 \boldsymbol{(}f,\,g\boldsymbol{)}{\lsm}_A \equiv \boldsymbol{(}f,\,g\boldsymbol{)}_{\!A(\vec{X})} =  \boldsymbol{(}\,\boldsymbol{(}K{(\vec{P},\vec{X})},\,f\boldsymbol{)}_{B(\vec{P})}\,,\,g\boldsymbol{)}_{\!A(\vec{X})} = \boldsymbol{(}\,\boldsymbol{(}K{(\vec{P},\vec{X})},\,g\boldsymbol{)}_{\!A(\vec{X})}\,,\,f\boldsymbol{)}_{\!B(\vec{P})}\\ = \boldsymbol{(}g,\,f\boldsymbol{)}_{B(\vec{P})} \equiv \boldsymbol{(}f,\,g\boldsymbol{)}{\ls}_{B} 
\end{multline}
(where the step in the middle corresponds to rearranging the integrals of the associated inner product expressions),
which proves the desired result for weighted norms.  This line of argument can be immediately generalized to include inner products with intrinsic differential operators, such as the gradient terms occurring in the (weighted) Dirichlet integral.

 Property (4) is usually proved in terms of Hilbert space inner products that have a quadratic form \cite{Aronszajn,BergmanSpaces2}, but the derivation just given in terms of inner products with an integral form is more natural in the present context since it clearly generalizes all the inner products encountered in \cite{DIDACKS} and it does not explicitly require the assumption of a countable basis.  (Other standard features of RKHS theory were also derived in \cite{DIDACKSII} without resorting to the usual technical functional analysis assumptions.)  Property~(4) means that there are usually infinitely many possible symmetric reproducing kernels for a given region since there a like number of possible norms.  Finally, routinely some norm and associated symmetric reproducing kernel  can be transformed by using the action of a positive definite linear differential operator (or linear representer) to yield a new norm and associated symmetric reproducing kernel (see, for example, Moritz \cite{Moritz}).

%For the following see Gratzer p. 253 \& ff., Hahn (LaTex for Everyone) p.199, etc.
%Here the number in braces indicates the size of the number field (1-9) => {9}; (10-99) => {99}.

\end{document}